\documentclass{amsart}
\usepackage{graphicx, amsmath, amssymb, amsthm}

\newtheorem{theorem}{Theorem}[section]

\theoremstyle{definition}

\theoremstyle{remark}
\newtheorem{remark}[theorem]{Remark}

\begin{document}

\title[It is not a Coincidence! On Curious Patterns in Optimization Problems]{It is not a Coincidence! \\ On Curious Patterns in 
Calculus Optimization Problems}
\author{Maria Nogin}

%\begin{biog}
%\item[Maria Nogin] (mnogin@csufresno.edu) received her PhD in mathematics from the 
%University of Rochester in 2003. She has been a faculty member at the California State University, Fresno since then. She is regularly teaching 
%Calculus and Problem Solving courses, and is fascinated by patterns and similarities occuring in varius problems. 
%\end{biog}

\maketitle

\begin{abstract}
In the first semester calculus course we learn how to solve optimization problems such as maximizing the volume of a box or a can 
given its surface area, i.e. the amount of material, or minimizing the surface area given the desired volume. Whether you are a 
student or a teacher, have you ever wished you knew the answer to a problem after one glance at it, without doing long calculations? The patterns 
shown and explained below will enable you to do it for a few ``standard'' problems, but more excitingly, they will illustrate some very beautiful 
problem solving techniques and give insights into how different areas of mathematics are connected.   
\end{abstract} 

\section{The basics} 

\subsection{Optimizing a rectangle.} 

It can be shown using the standard calculus optimization technique that of all rectangles with a given perimeter, the square has the largest 
area. (It also follows that of all rectangles with a given area, the square has the smallest perimeter.) There is also an easy geometric 
way to prove this. The picture below shows two rectangles: one of them, $s \times s$, is a square, and the other, $(s+h)\times(s-h)$, is not. They 
have the same perimeter ($4s$), but the square has a larger area since the areas of the shaded parts are equal, but the square also contains a 
region that the other rectangle does not.  

\begin{center}
\includegraphics[scale=0.5]{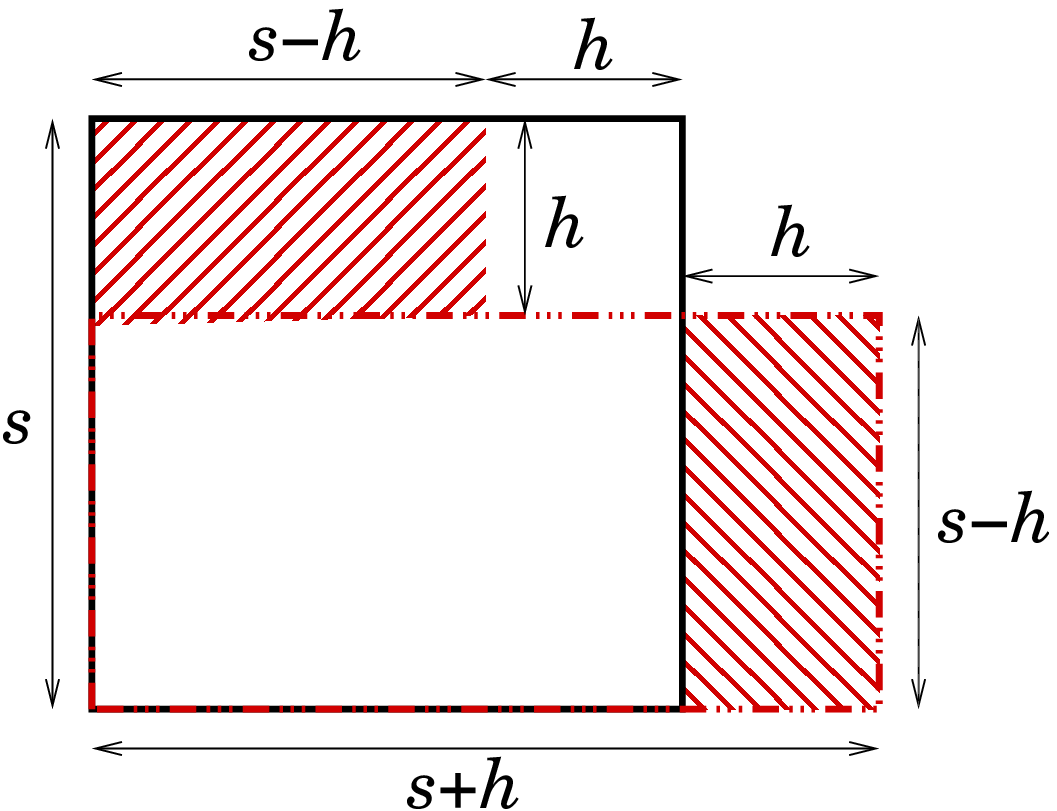}
\end{center}

\begin{remark} The above area inequality can also be seen algebraically: $$A_{square}=s^2>s^2-h^2=(s+h)(s-h)=A_{rectangle}$$ 
\end{remark} 

\begin{remark} As proved on p. 50 in \cite{niven}, even among all quadrilaterals with a given perimeter, the square has the largest area. 
\end{remark} 

\subsection{Optimizing a rectangular prism.} 

We can ask a similar question in 3D: of all rectangular prisms with the given volume, which one has the smallest surface area? Or, 
equivalently, of all rectangular prisms with the given surface area, which one has the largest volume? The answer to both of these questions is a 
cube. Notice that these problems are too general for a single variable calculus course: there are too many degrees of freedom. Namely, there are 
three unknowns (three dimensions), but only one constraint (the given surface area or the given volume). So typically in the first semester calculus 
class we only consider the 2-variable problem of optimizing a square prism. This would actually be sufficient for our work below, however, we can't 
resist the temptation to show our reader the following geometric argument that works for all rectangular prisms. Suppose we have any rectangular 
prism that is not a cube. Choose any face that is not a square and make it the base. Consider the square that has, say, the same area. It will have 
a smaller perimeter. Construct a prism with the same height as the original one:    

\begin{center}
\includegraphics[scale=0.6]{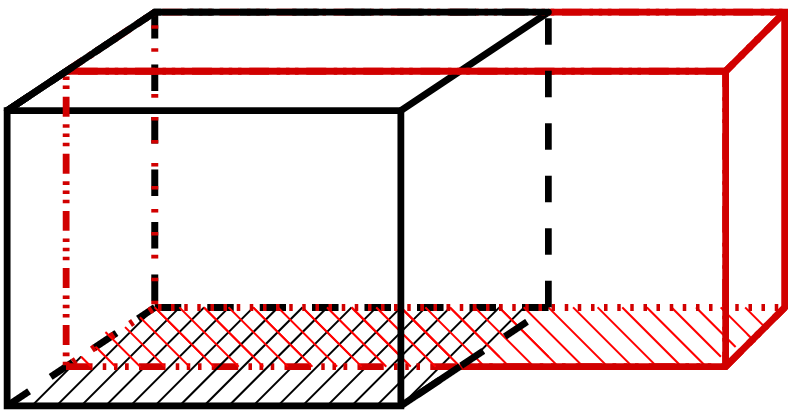}
\end{center}

It will have the same volume but a smaller surface area! (It will have the same volume because both the area of its base and its height are the same 
as those of the original prism; and a smaller surface area because the areas of the base and the top are the same but the area of the sides is 
smaller since the perimeter of the base is smaller.) Thus a prism that is not a cube cannot have the smallest possible surface area among all 
rectangular prisms with a fixed volume. 

Now we are ready to consider some more interesting problems. 

\section{The fence problem}

\subsection{Problem.} 

Some version of the following rectangular field optimization problem can be found in almost every calculus textbook (see e.g. 
exercise 7 on p. 337 in \cite{stewart}; example 1 on p. 258 in \cite{briggs} looks different at first, but actually is an equivalent 
optimization problem). 

\medskip 

\noindent\begin{tabular}{@{}p{8.6cm}@{\hspace{5mm}}p{3cm}@{}}
A farmer wants to fence off a rectangular field and divide it into 3 pens with fence parallel to one pair of sides. He has a total of 2400 ft of 
fencing. What are the dimensions of the field that has the largest possible area? & 
\raisebox{-1.8cm}{\includegraphics[scale=0.5]{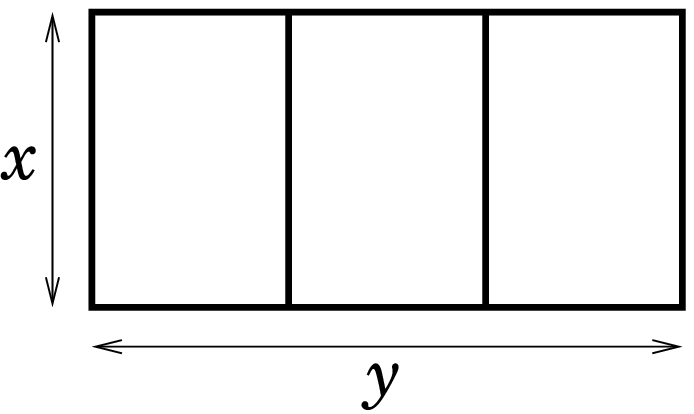}}
\end{tabular} 

The standard calculus optimization technique easily produces the answer: the field must have dimensions $x=300$ ft and $y=600$ ft. (It should be 
noted, however, that in this problem the area of the field turns out to be a quadratic function of $x$, therefore its maximum value, i.e. the value 
at the vertex of the parabola, can be found by completing the square, so this problem can be solved without using calculus.) 

\subsection{Observation.}

Notice that there are four pieces of fence of length 300 ft and two pieces of length 600 ft, therefore both the total length of the vertical 
pieces and the total length of the horizontal pieces are 1200 ft. The reader is invited to try a modification of the above problem (e.g. by varying 
the total length of the fence and/or the number of pens) and observe that the property described above still holds: to obtain 
the largest possible area, exactly half of the fence must be spent on the vertical pieces, and exactly half on the horizontal pieces. We could even 
have both vertical and horizontal partitions, and the optimal dimensions would still have the above property.    

\subsection{Functional explanation.} 
 
Let $L$ be the total length of the vertical pieces. Then $2400-L$ is the total length of the horizontal pieces. Then both dimensions are linear 
functions of $L$, namely, for the problem above, $x=\frac{L}{4}$ and $y=\frac{2400-L}{2}$. So the area is a quadratic function of $L$, namely, 
$\mbox{Area} = \frac{L}{4} \cdot \frac{2400-L}{2}$, but what is important, is that its roots are $L=0$ (well, if no fence at all is spent on 
the vertical pieces, the area of the field is 0) and $L=2400$ (similarly, if no fence is spent on the horizontal pieces, the area is 0).  
The graph of this quadratic function is a parabola opening downward, and, due to symmetry, the highest point on the parabola, its vertex, is halfway 
between the roots.  

\medskip 

\begin{center}
\includegraphics[scale=0.5]{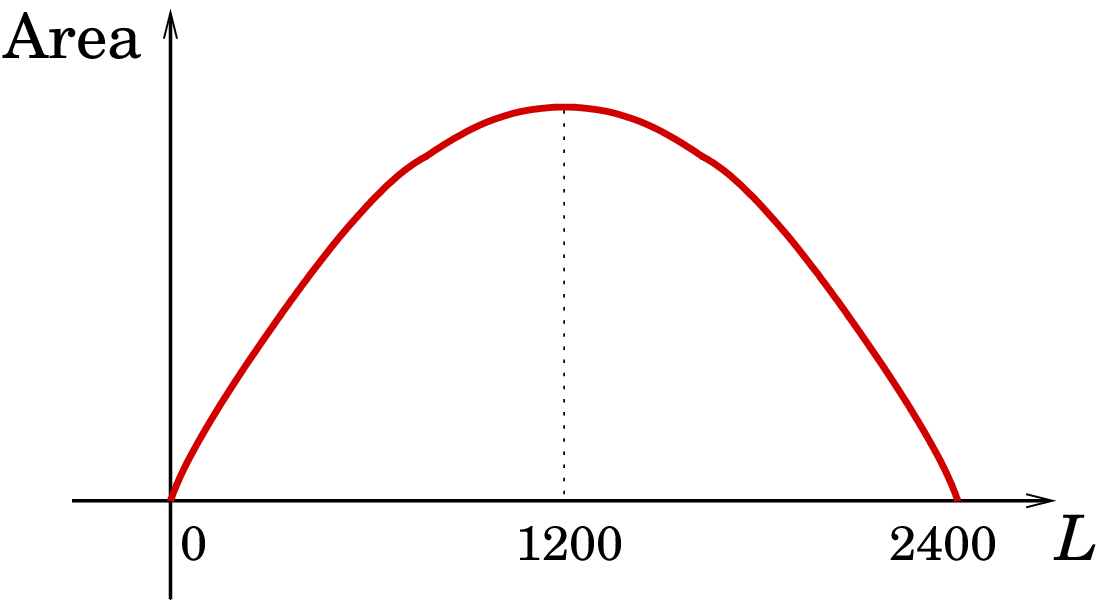} 
\end{center}

\subsection{Geometric explanation.} 

Let us consider the following modification of the problem: instead of having two partitions, we'll construct a rectangle that is twice ``taller'': 

\medskip 

\begin{center}
\includegraphics[scale=0.5]{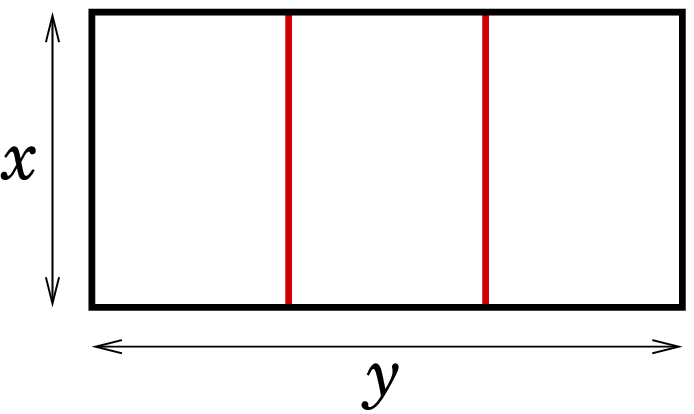} \hspace{1cm} \raisebox{1.2cm}{$\rightarrow$} \hspace{1cm} 
\includegraphics[scale=0.5]{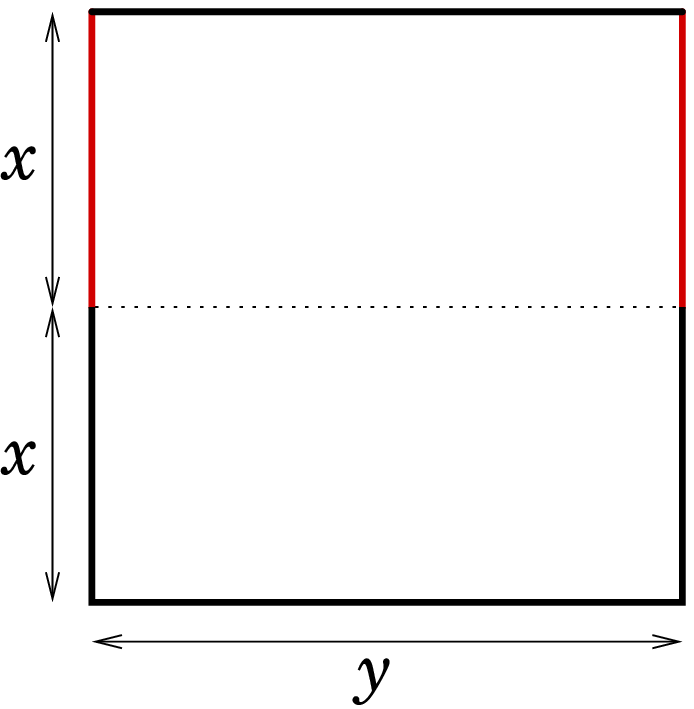}
\end{center}

This new rectangle uses exactly the same amount of fence as the original one, and its area is twice that of the original rectangle. Since maximizing 
the area is equivalent to maximizing twice the area, finding the optimal dimensions in the modified problem will give us the optimal dimensions in 
the original problem. But the modified problem is just a rectangle, and we know that the optimal shape is a square. In a square, the total length of 
the vertical pieces is equal to that of the horizontal pieces. Thus this property must hold in the original problem. 

We can use the same idea with any number of required partitions. Say, if we need four partitions to make five pens, we'll modify the problem by 
making the rectangle three times taller; from 6 partitions we'll make the rectangle four times taller, and so on. The reader might ask: what are we 
going to do if the number of partitions is odd? Well, we will just cut one of them in half: 
 
\medskip 

\begin{center}
\includegraphics[scale=0.5]{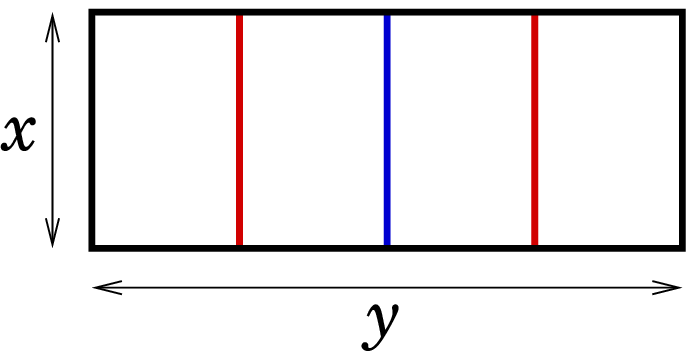} \hspace{1cm} \raisebox{1.2cm}{$\rightarrow$} \hspace{1cm} 
\includegraphics[scale=0.5]{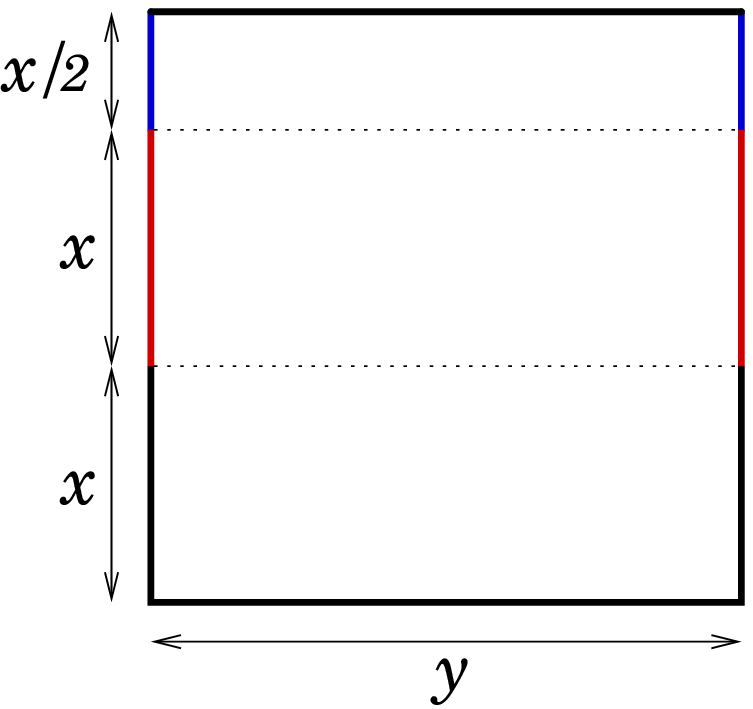}
\end{center}

No matter what the number of partitions is, the optimal solution is to spend exactly half of the fence on the vertical pieces and the other half on 
the horizontal pieces. 

\section{The cylindrical can problem} 

\subsection{Problem.} 

The following is another well-known problem (see e.g. example 2 on p. 333 in \cite{stewart}; problem 6 on p. 135 in \cite{silverman} is a 
slightly different version). 

\noindent\begin{tabular}{@{}p{9.7cm}@{\hspace{5mm}}p{3cm}@{}}
A cylindrical can must have volume 1000 cm$^3$. Find the dimensions that will minimize the cost of the metal to manufacture the can. In 
mathematical terms, find the dimensions that minimize the surface area of such a cylinder (with top and bottom). & 
\raisebox{-1.8cm}{\includegraphics[scale=0.5]{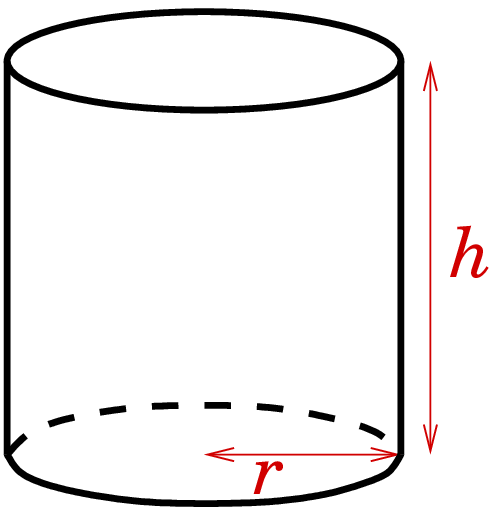}}
\end{tabular}

The reader is invited to check that the optimal dimensions are 
$r=\sqrt[3]{500/\pi}$ cm and $h=2\sqrt[3]{500/\pi}$ cm.  

\noindent\begin{tabular}{@{}p{8.9cm}@{\hspace{5mm}}p{3cm}@{}} 
\subsection{Observation.} Notice that in the optimal can $h=2r$, and the same relationship holds for an optimal square prism box if we let $r$ be 
the distance from the center of the base to any side (note that $r$ is the radius of the inscribed circle, called the inradius of the square; the 
line segment from the center of a regular polygon to the middle of any side is also called its apothem). & 
\raisebox{-2.8cm}{\includegraphics[scale=0.5]{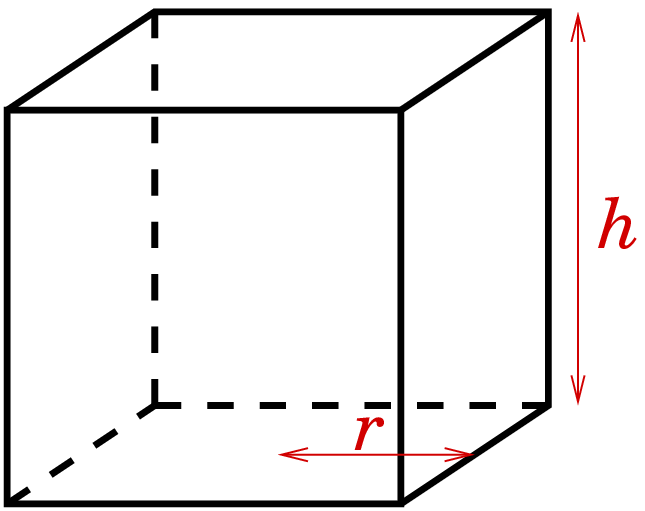}} 
\end{tabular}  

Is this a coincidence? We will show that it is not, but rather these two problems (optimizing the can and optimizing the square prism) are 
equivalent. 

\subsection{Why is this so?} 

Let us see how the volume and the surface area of the can and the square prism are related: 
$$V_{can} = A_{circle} h, \ \ \ \ \ \ V_{prism} = A_{square} h.$$ 
So $$V_{prism}=\frac{A_{square}}{A_{circle}} V_{can} = \frac{4r^2}{\pi r^2} V_{can} = \frac{4}{\pi} V_{can}.$$ 
Now, 
$$SA_{can} = 2 A_{circle} + P_{circle}h, \ \ \ \ \ \  SA_{prism} = 2 A_{square} + P_{square}h.$$

Notice that the expressions for both surface areas contain twice the area of the base, and 
$A_{square} = \frac{4}{\pi} A_{circle}$. 
If only the other two terms in these expressions were related in the same way, namely, if $P_{square}$ was $\frac{4}{\pi}$ times larger than 
$P_{circle}$, that would mean that $SA_{prism}$ was also $\frac{4}{\pi}$ times larger than $SA_{can}$, which, in turn, would mean that the two 
problems are equivalent. 

So here is our question: is $\frac{P_{square}}{P_{circle}} = \frac{A_{square}}{A_{circle}}$?
It is easy to see that these quotients are $\frac{8r}{2\pi r}$ and $\frac{4r^2}{\pi r^2}$, so indeed they are equal. 
But is that a coincidence, or is there something deeper here? Why $\frac{P_{square}}{P_{circle}} = \frac{A_{square}}{A_{circle}}$?

Let us rewrite the above equation as follows: $\frac{A_{circle}}{P_{circle}} = \frac{A_{square}}{P_{square}}$.  

Now look again at the above equation in terms of $r$: $\frac{\pi r^2}{2\pi r} = \frac{4r^2}{8r}$. On each side, the denominator is 
the derivative of the numerator! This property is briefly illustrated below, while its formal proof for any regular polygon is given in 
\cite{miller} and a more general case is explored in \cite{dorff}. 

Imagine increasing the radius of a circle by $\Delta r$ (see 
picture below). Its area is increased by $\Delta A \approx P_{circle} \Delta r$. That's why the derivative of the area is $\displaystyle 
\lim_{\Delta r\to 0}\frac{\Delta A}{\Delta r} = \lim_{\Delta r\to 0} \frac{P_{circle} \Delta r}{\Delta r} = P_{circle}$. The same is true for the 
square! 

\medskip 

\begin{center}
\includegraphics[scale=0.35]{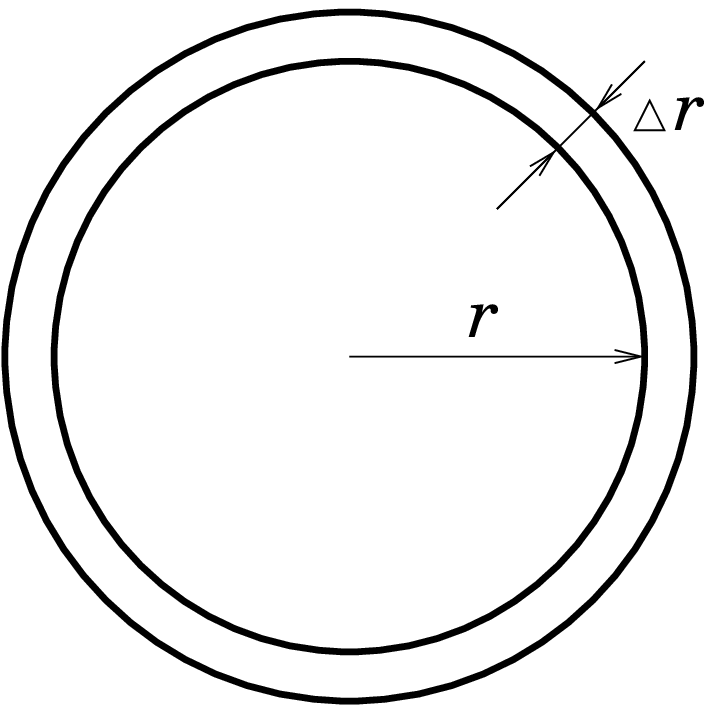} \hspace{1cm} \includegraphics[scale=0.35]{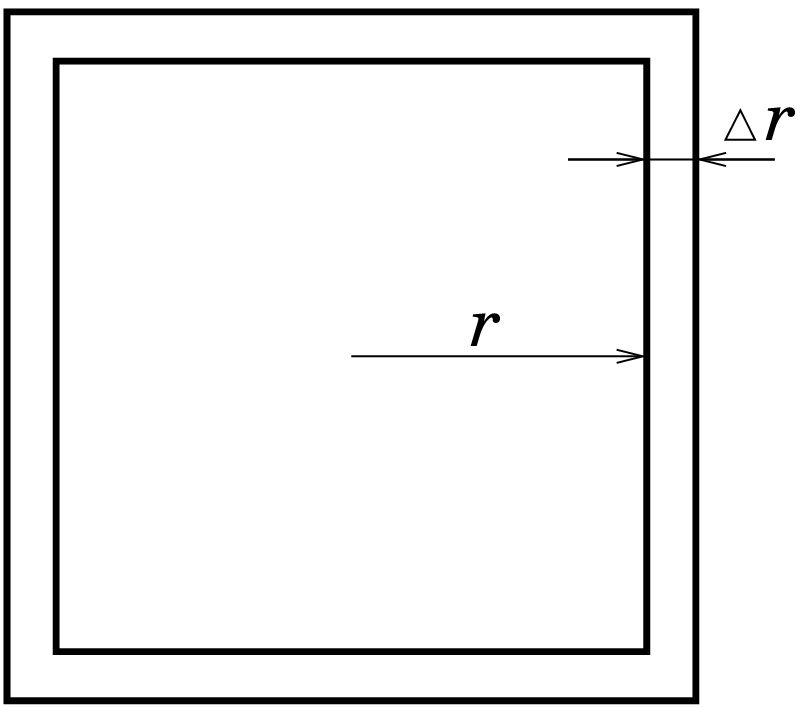}
\end{center}

Now, both areas (of the circle and of the square) are quadratic power functions of $r$. For any two quadratic power functions, say, $ar^2$ and 
$br^2$, their derivatives are $2ar$ and $2br$ respectively, so the ratios of the function to its derivative would be 
$\frac{ar^2}{2ar}=\frac{r}{2}=\frac{br^2}{2br}$. It does not matter what the coefficient $a$ or $b$ is, since it appears in both the function and 
its derivative and gets canceled in the ratio.    

So the above observation is not at all a coincidence. Moreover, as we will see in the next section, this applies to any prism that has a 
regular polygon as the base! 

\subsection{Other boxes.} 

Suppose we want to construct a box in the shape of some other prism whose base is a regular polygon. What would be the optimal shape then? For 
example, given the desired volume, what ratio of the height $h$ to the inradius $r$ would minimize the surface area (and thus minimize the amount of 
material needed)?  

\begin{center}
\includegraphics[scale=0.45]{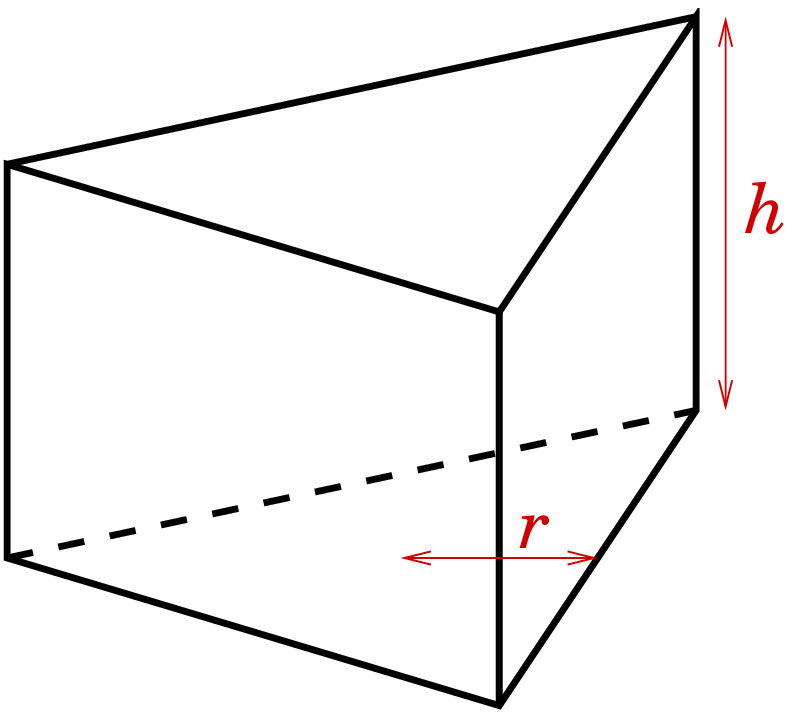} \hspace{1.5cm} \raisebox{0.3cm}{\includegraphics[scale=0.5]{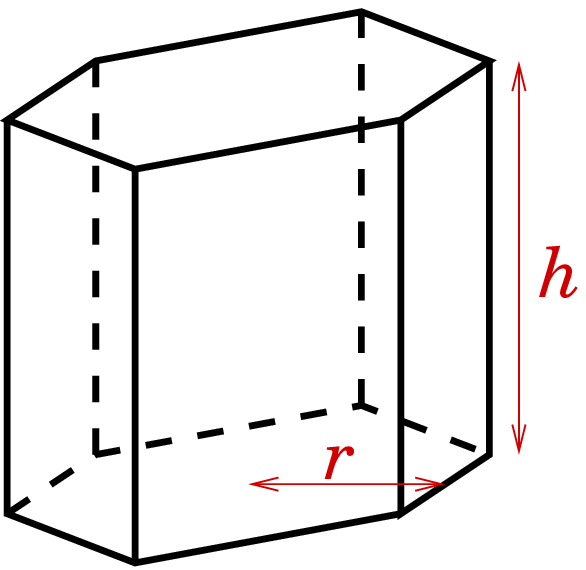}} 
\end{center} 

As we know, area grows proportionally to the square of the length. Therefore, the area of any regular polygon is a square function of the 
inradius, i.e. of the form $cr^2$ where $c$ is some constant. 
The perimeter of the polygon is then $2cr$. This means that both the volume and the surface area of any such prism differ from those of a square 
prism by a constant multiple (namely, $\frac{c}{4}$), thus the problem of optimizing any such box is 
equivalent to the problem of optimizing a square prism. So, in all cases the optimal shape is the one for which $h=2r$!  

\section{The ellipse inscribed in a semi-circle problem} 

The following is a modification of problem 9 on p. 979 in \cite{stewart}. 

\vspace{-3mm}

\noindent\begin{tabular}{@{}p{7.1cm}@{\hspace{5mm}}p{5cm}@{}}
\subsection{Problem.} 
Of all ellipses inscribed in a semi-circle of radius 1, find the one with the largest possible area.  \phantom{mmmmmmmmmmm}
Hint: If the semicircle is given by the equation $x^2+y^2=1$, $y\ge 0$, the ellipse should have equation of the form 
$\frac{x^2}{a^2}+\frac{(y-b)^2}{b^2}=1$.& 
\raisebox{-3.3cm}{\includegraphics[scale=0.5]{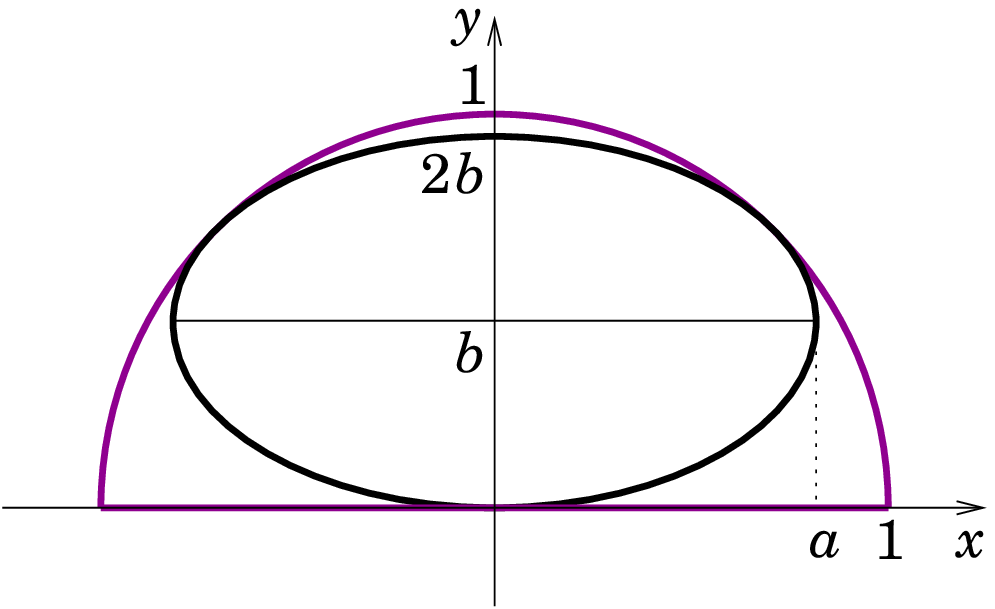}} 
\end{tabular}  

This problem is a bit tougher than the previous ones, and requires a longer computation. The optimal shape turns out to be the one with
$a=\frac{\sqrt{6}}{3}$ and $b=\frac{\sqrt{2}}{3}$.

\subsection{Observation.}

Let us dig a bit deeper and calculate the points of intersection of the ellipse with the semicircle. They are $\left( \pm \frac{\sqrt{2}}{2}, 
\frac{\sqrt{2}}{2} \right)$. Notice that these are exactly vertices of the largest-area rectangle inscribed in a semi-circle (see e.g. example 5 on 
p. 336 in \cite{stewart}; see also p. 130 in \cite{niven} for a generalization of this problem):   

\begin{center}
\includegraphics[scale=0.5]{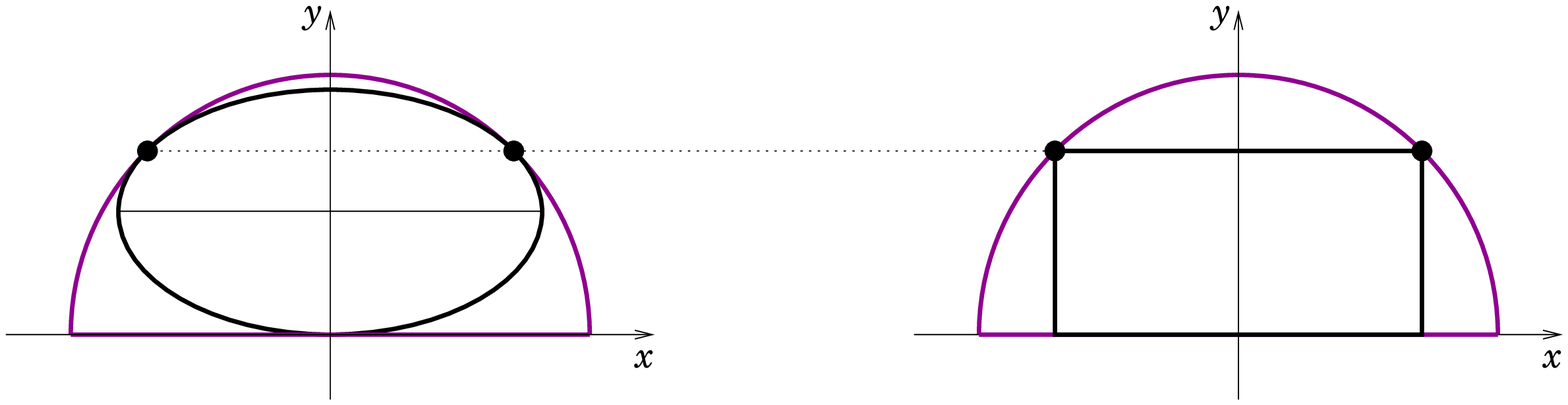}
\end{center}

\subsection{Why	is this	so?} Are the problems of optimizing the ellipse and optimizing the rectangle inscribed in a semi-circle equivalent? It 
certainly appears so, but what do the inscribed ellipse and the inscribed rectangle have in common? Try to figure it out! 

{\bf Acknowledgments.}  
The author would like to thank her colleagues at Fresno State University for numerous insightful conversations about Calculus problems as well as 
for their encouragement to write this paper, and the anonymous reviewer for many good suggestions on its improvement.   

%\begin{abstract}
%In this paper we consider a few Calculus optimization problems in which we notice peculiar patterns. In each of these cases there is a 
%geometric explanation for the pattern showing that it is not just a coincidence.
%\end{abstract}

\vfill\eject

\end{document}